\documentclass[12pt,a4paper]{article}

\usepackage{theorem}
\usepackage{enumerate}
\usepackage{amsmath}
\usepackage{latexsym}
\usepackage{amssymb}
\usepackage{amsfonts}
\usepackage{graphicx, epsfig,epstopdf, subfigure}
\usepackage[varg]{txfonts}
\usepackage{here}
\usepackage{multirow}
\usepackage{array}

\newcounter{RomanNumber}
\newcommand{\MyRoman}[1]
{\setcounter{RomanNumber}{#1}\Roman{RomanNumber}}

\theorembodyfont{\normalfont\itshape}
\newtheorem{thm}{Theorem}

\newtheorem{Thm}{Theorem}

\newtheorem{lem}{Lemma}

\newtheorem{clm}{Claim}[section]

\newtheorem{fact}{Fact}
\newtheorem{conj}{Conjecture}[section]

\theoremstyle{remark}

{\theorembodyfont{\upshape}
}

\newcommand{\proof}{\medbreak\noindent\textit{\textbf{Proof.}}\quad}

\newcommand{\qed}{{$\quad\square$\vspace{3.6mm}}}

\bfseries\normalfont

\begin{document}

\title{On directed version of the Sauer-Spender Theorem \thanks{The author’s work is supported by NNSF of China (No.11671232).}}

\author{Yun Wang, Jin Yan\thanks{Corresponding author: Jin Yan, Email:  yanj@sdu.edu.cn} \unskip\\[.5em]
School of Mathematics, Shandong University, Jinan 250100, China}

\date{}
\maketitle

\vspace{-12pt}
\begin{abstract}
Let $D=(V,A)$ be a digraph of order $n$ and let $W$ be any subset of $V$. We define the minimum semi-degree of $W$ in $D$ to be $\delta^0(W)=\mbox{min}\{\delta^+(W),\delta^-(W)\}$, where $\delta^+(W)$ is the minimum out-degree of $W$ in $D$ and $\delta^-(W)$ is the  minimum in-degree of $W$ in $D$. Let $k$ be an integer with $k\geq 1$. In this paper, we prove that for any positive integer partition $|W|=\sum_{i=1}^{k}n_i$ with $n_i\geq 2$ for each $i$, if $\delta^0(W)\geq \frac{3n-3}{4}$, then there are $k$ vertex disjoint cycles $C_1,\ldots,C_k$ in $D$ such that each $C_i$ contains exactly $n_i$ vertices of $W$. Moreover, the lower bound of $\delta^0(W)$ can be improved to $\frac{n}{2}$ if $k=1$, and $\frac{n}{2}+|W|-1$ if $n\geq 2|W|$.

The minimum semi-degree condition $\delta^0(W)\geq \frac{3n-3}{4}$ is sharp in some sense and this result partially  confirms the conjecture posed by Wang [Graphs and Combinatorics 16 (2000) 453-462]. It is also a directed version of the Sauer-Spender Theorem on vertex disjoint cycles in graphs [J. Combin. Theory B, 25 (1978) 295-302].

\end{abstract}

\vspace{-6mm}
\section{Introduction}  \label{section 1}
In this paper, we consider finite simple graphs or digraphs, which have neither loops nor multiple
edges or arcs.  For terminology and notation not defined in this paper, we refer the readers to \cite{Bang2009} and \cite{Bondy1976}.  For a graph $G=(V,E)$ (or a digraph $D=(V,A)$), we denote by $V$ the vertex set of $G$ (or $D$) and the cardinality of $V$ is the \emph{order} of $G$ (or $D$). The edge set of $G$ or the arc set of $D$ is defined by $E$ or $A$, respectively.

For a vertex $v$ in $D$, let the degree of $v$  to be $d_D(v)=d^+_D(v)+d^-_D(v)$. The minimum degree of $D$ is often written by $\delta(D)$ and the \emph{minimum semi-degree} of $D$  is $$\delta^0(D)=\mbox{ min }\{\delta^+(D),\delta^-(D)\},$$  where $\delta^+(D)$, $\delta^-(D)$ is the minimum out-degree, minimum in-degree of $D$, respectively.

For a set $W\subseteq V$, the \emph{minimum degree of $W$ in $G$}, is defined by $\delta(W)=\mbox{ min }\{d_G(v):v\in W\}$. The \emph{minimum out-degree of $W$ in $D$} is  $\delta^+(W)=\mbox{ min }\{d^+_D(v):v\in W\}$. Similarly, one can define $\delta^-(W)$, and the \emph{minimum semi-degree of $W$ in $D$} is $\delta^0(W)=\mbox{ min }\{\delta^+(W),\delta^-(W)\}$. A cycle (path) in digraphs means a directed cycle (path). We say $k$ pairwise vertex disjoint cycles by simply saying $k$ disjoint cycles unless otherwise specified. We denote by $G^*$ the symmetric digraph obtained from $G$ by replacing every edge $xy$ with the pair $xy,yx$ of arcs.

Let $k$ be an integer with $k\geq 1$. For every integer partition $|W|=n_1+\cdots+n_k$ with $n_i\geq 3$ for each $i$, if $G$ contains $C_{1}\cup\cdots\cup C_{k}$ as a subgraph with $|V(C_i)\cap W|=n_i$ ($1\leq i\leq k$), then we say that $G$ has \emph{an arbitrary $W$-cycle-factor.}  We abbreviate $W$-cycle-factor to \emph{cycle-factor} or \emph{2-factor} when $W=V(G)$. In the same way, we can define \emph{an arbitrary $W$-cycle-factor of a digraph $D$}, where $n_i\geq 2$.

\medskip

In 1952, Dirac \cite{Dirac1952} proved that every graph $G$ of order $n\geq 3$ with $\delta(G)\geq n/2$ is hamiltonian,  which is well known in graph theory.  As a generalization of this,  Corr\'{a}di and Hajnal \cite{Corradi1963} pointed out that a graph of order $n\geq 3k$ with $\delta(G)\geq2k$ contains $k$ disjoint cycles. 

When the cycles are required to be of the specific lengths,  the problem becomes much complicated. Sauer and Spender \cite{Sauer1978} first gave the following degree condition for graphs to contain an arbitrary 2-factor.

\begin{Thm}\label{sauerthm} \cite{Sauer1978} Let $G$  be a graph of order $n$. If $\delta(G)\geq 3n/4$, then $G$ contains an arbitrary 2-factor.\end{Thm}

They further conjectured that the same conclusion can be guaranteed by minimum degree at least $2n/3$ and, it was verified by Aigner and Brandt.

\begin{Thm}\label{thm3} (\cite{Aigner1993}) A graph $G$  of order $n$ with $\delta(G)\geq (2n-1)/3$  contains an arbitrary 2-factor.\end{Thm}

Let $W$ be a set of $V$. It is called \emph{cyclable} if $G$ or $D$ has a $W$-cycle-factor with exactly one cycle. Bollob\'{a}s and Brightwell \cite{Bollobas1993} showed that the set $W$ of $V(G)$ with  $\delta(W)\geq n/2$ is cyclable, which  extended  the Dirac's classical theorem. In \cite{Hong2015},  Wang introduced a new way to generalize Theorem \ref{thm3}  and conjectured as follows.
\begin{conj}\label{conj5}(\cite{Hong2015}) Let $k$ be an integer.  Suppose that $G$ is a graph of order $n$ and $W$ is a set of $V(G)$ with $|W|\geq 3k$. If $\delta(W)\geq 2n/3$, then $G$ contains an arbitrary $W$-cycle-factor.\end{conj}

In the same paper, he proved that under the same condition as in Conjecture \ref{conj5} the graph $G$ contains $k$ disjoint cycles covering $W$, such that each of $k$ cycles contains at least three vertices of $W$.

Generally speaking, the problem of finding a cycle through a given vertex set is more difficult for digraphs than for  graphs, even for the special case of two vertices. For example, Fortune, Hopcroft and Wyllie \cite {Fortune1980} have shown that the problem of finding a cycle through two prescribed vertices in digraphs is $NP$-complete, while the same problem for undirected graphs is relatively easy.

In \cite{Houri1962}, Ghouila-Houri gave the directed version of  Dirac's result as follows: A digraph $D$ of order $n$ with $\delta^0(D)\geq n/2$ is hamiltonian. Berman and  Liu \cite{Berman1998} considered a minimum degree sum condition for strong digraphs to be cyclable.  
In this paper, we consider the minimum semi-degree condition in digraphs. More precisely, we prove the following.

\begin{thm}\label{thmcyclable} Let $D$ be a digraph of order $n$ and let $W\subseteq V(D)$. If $\delta^0(W)\geq n/2$, then $W$ is cyclable.\end{thm}

The complete bipartite digraph $K_{k,k+1}^*$ shows that the degree condition in Theorem \ref{thmcyclable} is best possible.

For digraphs to contain $k$ disjoint cycles, Wang \cite{Hong2000}  raised the following two conjectures, where the second one extended Theorem \ref{thm3} to digraphs.

\begin{conj}  (\cite{Hong2000}) \label{conj4} Let $k$ be an integer and let $D$ be a digraph of order $n\geq 3k$. If $\delta^0(D)\geq (9k-3)/4$, then $D$ contains $k$ disjoint cycles of order at least $3$.\end{conj}

\begin{conj} (\cite{Hong2000}) \label{conj3} A digraph $D$  of order $n$ with $\delta(D)\geq(3n-3)/2$ contains an arbitrary cycle-factor. \end{conj}


Subsequently, Wang \cite{Hong2013} proposed  Conjecture \ref{conj1}, which is a generalization of
Conjectures \ref{conj5} and \ref{conj3}.


\begin{conj} (\cite{Hong2013}) \label{conj1}%
Let $D$ be a digraph of order $n\geq 6$. Suppose that $W$ is a set of $V(D)$ with $\delta(W)\geq (3n-3)/2$. Then $D$ contains an arbitrary $W$-cycle-factor.\end{conj}

Conjecture \ref{conj3} is essentially the case $W=V(G)$ of Conjecture \ref{conj1}. In \cite{Hong2013}, Wang claimed that the case $k=2$ holds.  This paper proves that it holds for minimum semi-degree condition by proving Theorem \ref{thm1}.
\vspace{-2mm}



\begin{thm}\label{thm1} Let $D$ be a digraph with $W\subseteq V(D)$. If  $\delta^0(W)\geq (3n-3)/4$, then $D$ contains an arbitrary $W$-cycle-factor.\end{thm}

Clearly, Theorem \ref{thm1} implies the case $n=3k$ of Conjecture \ref{conj4} and the case $\delta^0(D)\geq (3n-3)/4$ of Conjecture \ref{conj3}. The degree condition in  Theorem \ref{thm1} is best possible in some sense (see Remark in Section \ref{section 7}).  Here, we further present a result which optimizes the bound of the minimum semi-degree if $|V(D)|$ is sufficiently large compared with $|W|$.

\begin{thm}\label{thm5} If $n\geq 2|W|$, then the minimum semi-degree condition in Theorem \ref{thm1} can be replaced by $\delta^0(W)\geq \frac{n}{2}+|W|-1$.\end{thm}

Several other conjectures on the existence of disjoint cycles have appeared in the literature. Let us mention one of the most interesting, which is still open. El-Zahar \cite{Zahar1984} conjectured that if $G$ is a graph with $\delta(G)\geq\frac{n+C}{2}$, where $C$ is the number of odd $n_i$s, then $G$ has an arbitrary 2-factor. For results on this conjecture, see \cite{Abbasi1998,Corradi1963,Hong2010,Hong2012}.

\medskip
Let $W$ be a set of $V(G)$. Table \ref{table1} summarizes the results mentioned above, where "$k$ disjoint cycles" in the fourth column means that $G$ or $D$ has $k$ disjoint cycles in which each cycle contains at least three vertices of $W$.  Recall that if $W=V(G)$, then a $W$-cycle-factor is actually a cycle-factor.

\begin{table}[H]\label {table1}
\centering\large
\resizebox{15cm}{3.5cm}{
\begin{tabular}{|p{2cm}<{\centering}||p{3cm}<{\centering}||p{6cm}<{\centering}|p{6cm}<{\centering}|p{9cm}<{\centering}|}
\hline\hline
&{}&{$W$-cycle-factor with 1 cycle}&{$k$ disjoint cycles}&{$W$-cycle-factor with $k$ ($\geq 2$) cycles}\\
\hline\hline
\multirow{5}{1.2cm}{graph}  &\multirow{3}{2.5cm}{ degree condition }&\multirow{3}{3cm}{$\delta(G)\geq n/2$ (Dirac \cite{Dirac1952})} &\multirow{3}{5cm}{$n\geq 3k,\delta(G)\geq 2k$ (Corr\'{a}di-Hajnal \cite{Corradi1963})}&{$\delta(G)\geq 3n/4$ \,\,(Sauer-Spencer \cite{Sauer1978}}) \\
\cline{5-5}                 & & & &{$\delta(G)\geq (2n-1)/3$ \,\,(Aigner-Brandt \cite{Aigner1993}})\\
\cline{5-5}                 & & & &{$\delta(G)\geq(n+C)/2$ \,\,(El-Zahar \textbf{Conj.} \cite{Zahar1984})}\\

\cline{2-5}                 &\multirow{2}{2cm}{ partial degree condition}&\multirow{2}{6cm}{$\delta(W)\geq n/2$ (Bollob\'{a}s-Brightwell \cite{Bollobas1993})}&\multirow{2}{3cm}{$\delta(W)\geq 2n/3$ (Wang \cite{Hong2015})}                &\multirow{2}{5cm}{$\delta(W)\geq 2n/3$ \,\,\,\,\,\,\,\,\,\,\,\,\,\,\,\,\,(Wang \textbf{Conj.} \cite{Hong2015})}\\
& & &  &\\
\hline
\multirow{5}{1.2cm}{digraph}  &\multirow{2}{2.5cm}{ degree condition }&\multirow{2}{5cm}{$\delta^0(D)\geq n/2$ (Ghouila-Houri \cite{Houri1962})} &\multirow{2}{6cm}{$n\geq 3k,\delta^0(D)\geq (9k-3)/4$ \,\,\,\,\,\,\,\,\,(Wang \textbf{Conj.} \cite{Hong2000})}&\multirow{2}{5cm}{$\delta(D)\geq (3n-3)/2$ (Wang \textbf{Conj.} \cite{Hong2000})} \\
     & & &
& \\
\cline{2-5}                 &\multirow{3}{2cm}{ partial degree condition}&\multirow{3}{4cm}{$\delta^0(W)\geq n/2$\\(\textbf{\textit{Theorem \ref{thmcyclable}}})}& &{$\delta(W)\geq (3n-3)/2$(Wang \textbf{Conj.} \cite{Hong2013})}\\
\cline{5-5}
& &
 &  &{$\delta^0(W)\geq (3n-3)/4$} (\textbf{\textit{Theorem \ref{thm1}}})\\

 \cline{5-5}
& &
 &  &{$n\geq 2|W|$,\,$\delta^0(W)\geq n/2+|W|-1$}(\textbf{\textit{Theorem \ref{thm5}}})\\
\hline
\end{tabular}}\caption{Comparisons of results.}
\end{table}

The rest of the paper is organized as follows: The main aim of Section \ref{section 2} is to establish some notation and,  present the equivalent form (Theorem \ref{thm2}) of Theorem \ref{thm1} in bipartite graphs. In Section \ref{section 4}, we show some lemmas which are useful in the proof of Theorem \ref{thm2}.  In order to maintain the consistency of the article, we prove Theorem \ref{thm2} first  in Section~\ref{section 5}.  The proofs of  Theorems \ref{thmcyclable} and \ref{thm5} are presented in Section~\ref{section 6}. Section~\ref{section 7} discusses the degree condition in Theorem \ref{thm1}.
\vspace{-6mm}

%
\section{Preparation}  \label{section 2}

We prepare terminology and notation which will be used in our proofs. For $X \subseteq V$, we use $[X]$ to denote the subgraph (subdigraph) induced by $X$.  For disjoint subsets or subgraphs $X$ and $Y$ of $G$, define $e(X,Y)$  to be the number of edges between $X$ and $Y$.

We use $M$ to denote a perfect matching of $G$.  An $M$-\emph{alternating cycle (path)} is a cycle (path) such that edges belong to $M$ and not to $M$, alternatively. Let $M^{\prime}$ be a  subset of $M$.  For convenience, we use $M^{\prime}(H)$  to denote $M^{\prime}\cap E(H)$ for any subgraph or subset $H$ of $G$ and the matching $M^{\prime}$. The $M^{\prime}$-\emph{length of} $H$, i.e. the number of edges in $M^{\prime}(H)$, is defined by $|M^{\prime}(H)|$.

Let  $C=v_1v_2\cdots v_{l}v_1$ be a cycle in $G$ or $D$. The \emph{length} of a cycle $C$, i.e. the number of edges of $C$, is written by $l(C)$. The successor of $v_i$ on $C$, written as $v_i^+$, is the vertex $v_{i+1}$. Similarly, one can denote the predecessor of $v_i$ by $v_i^-$. We denote $C[v_i,v_j]=v_iv_{i+1}\cdots v_j$. (The subscripts are taken modulo in $\{1,2,\ldots,l\}$). It is trivial to extend the above definitions to a path.

\bigskip



To obtain short proofs of various results on digraph $D = (V, A)$, the following transformation to the class of bipartite graphs is extremely useful (in \cite{Bang2009}). Let $BG(D) = (V_X, V_Y; E)$ denote the bipartite graph with partite sets $V_X,V_Y$ such that $V_X= \{v_X : v \in V\}, V_Y= \{v_Y : v \in V\}$. The edge set $E$ of $BG(D)$ is defined by  $E=\{u_Xv_Y: u_X\in V_X, v_Y\in V_Y \mbox{ and } uv\in A\}$. We call $BG(D)$ the bipartite representation of $D$.

We construct a bipartite graph $G(X,Y)$ from $BG(D)$ by adding the perfect matching $M=\{v_Xv_Y: v\in V\}$. The method of transforming a digraph $D$ into $G(X,Y)$  is also in \cite{Chiba2018} and \cite{Zhang2013}. Note that, in this construction, the following properties hold:

$\bullet$ $uv\in A(D)$ if and only if $u_Xv_Y\in E(G(X,Y))$. Moreover, $d_G(v_X)=d^+_D(v)+1$ and $d_G(v_Y)=d^-_D(v)+1$, and

$\bullet$ an $M$-alternating cycle of length $2l$ $(\geq 4)$ in $G(X,Y)$ corresponds to a directed cycle of length $l$ $(\geq 2)$ in $D$.

\medskip
We can now rephrase Theorem \ref{thm1} in bipartite graphs as follows.

\begin{thm}\label{thm2} Let $G(X,Y)$ be a balanced bipartite graph with perfect matching $M=\{x_iy_i: x_i\in X, y_i\in Y, 1\leq i\leq n\}$. Suppose that $M_0\subseteq M$ and $\delta(M_0)\geq (3n+1)/4$. Then for every positive integer partition $|M_0|=\sum_{i=1}^kn_i$ with $n_i\geq 2$ $(1\leq i\leq k)$,  $G(X,Y)$ contains $k$ disjoint $M$-alternating cycles $C_1,\ldots, C_k$ such that each $C_i$ contains exactly $n_i$ edges of $M_0$.\end{thm}
\medskip

In the following, let $M$ be a perfect matching of $G(X,Y)$ and $M_0$  a subset of $M$. Set $M_1=M\backslash M_0$. We say an $M$-alternating cycle $C$ of $M_0$-length $l_{M_0}(C)\geq 2$ by simply saying a \emph{feasible cycle} unless otherwise specified. A path $P=x_1y_1\cdots x_ry_r$ is a \emph{feasible path} if it is an $M$-alternating path with $\{x_1,y_r\}\subseteq V(M_0)$. Moreover, for any feasible path $P^{\prime}=x^{\prime}_1y^{\prime}_1\cdots x^{\prime}_ry^{\prime}_r$ in $[V(P)]$ with $M_0(P)=M_0(P^{\prime})$, if $$e(\{x_1,y_r\},P)\leq e(\{x^{\prime}_1,y^{\prime}_r\},P^{\prime}),$$ then $P$ is a \emph{good feasible path}.
 Note that for any $M$-alternating path $P$, one can choose a good feasible path in $[V(P)]$.
\medskip

In order to simplify the proofs, we have the following  hypothesis:

\begin{fact}\label{fact1} For arbitrary two edges $f_i,f_j\in M_1(G)$, we may assume that there is no edge between $f_i$ and $f_j$.\end{fact}
Indeed, if there is an edge $f$ between $f_i$ and $f_j$, then we consider the graph $G$ which obtained by deleting the edge $f$ from the original graph. If $G$ has $k$ desired cycles, then  the same conclusion also holds in the original graph.
\vspace{-6mm}
\section{Lemmas}  \label{section 4}

In this section, the graph $G$ will refer to the bipartite graph in the Theorem 4 hypothesis unless mentioned otherwise.

\begin{lem}\label{lem1} If $C=x_1y_1\cdots x_ry_rx_1$ is a feasible cycle and $xy\in M_0\backslash E(C)$ with $e(\{x,y\},C)\geq r+1$, then  $[V(C)\cup\{x,y\}]$ contains a feasible cycle $C^{\prime}$ such that $M_0(C^{\prime})=M_0(C)\cup\{xy\}$.
\end{lem}
\proof  Set $M(C)=\{x_iy_i\,|\,1\leq i\leq r\}$. Suppose the lemma is false. Then $e(\{x,y\},\{y_i,x_{i+1}\})\leq 1$ for each pair $y_i,x_{i+1}$, where the subscripts $i$ are taken modulo in $\{1,\ldots, r\}$.  It follows that $e(\{x,y\},C)=\sum_{i=1}^r e(\{x,y\},\{y_i,x_{i+1}\})\leq r$, a contradiction.\qed

\begin{lem}\label{lem2} Let $s$ and $t$ be two integers such that $t\geq s\geq 2$ and $t\geq 3$. Suppose  $C_1$ and $C_2$ are two disjoint feasible cycles of $G(X,Y)$ such that $l_{M_0}(C_1)=s$, $l_{M_0}(C_2)=t$ and $$\sum_{f\in M_0(C_2)}e(f,C_1)>\frac{3}{4}tl(C_1).$$ Then $[V(C_1\cup C_2)]$ contains two disjoint feasible cycles $C^{\prime}$ and $C^{\prime\prime}$ such that $l(C^{\prime})+l(C^{\prime\prime})<l(C_1)+l(C_2)$.\end{lem}

\proof By way of contradiction, suppose the lemma fails. Let $C_1$ and $C_2$ be two cycles satisfying the condition of the lemma with $l(C_1)+l(C_2)$ minimum such that $[V(C_1\cup C_2)]$ does not contain the two required cycles $C^{\prime}$ and $C^{\prime\prime}$. We claim that
\begin{eqnarray}\label{eq1}
l(C_2)=2t, \,\, i.e. \,\, M_1(C_2)=\emptyset.
\end{eqnarray}

\textbf{\textit{Proof of (\ref{eq1}).}} Recall the definition of $v_i^+$ and $v_i^-$. Suppose (\ref{eq1}) does not hold and say $x_iy_i\in M_1(C_2)$, then we add a new edge $x_i^-y_i^+$ to $C_2-x_iy_i$ to obtain a feasible cycle $C$ as $t\geq 3$.  Clearly, $\sum_{f\in M_0(C)}e(f,C_1)>\frac{3}{4}tl(C_1)$. Then by the minimality of $l(C_1)+l(C_2)$, we derive that $[V(C_1\cup C)]$ contains two disjoint feasible cycles $L_1$ and $L_2$ with $l(L_1)+l(L_2)<l(C_1)+l(C)=l(C_1)+l(C_2)-2.$ Replacing $x_i^-y_i^+$ by $x_i^-x_iy_iy_i^+$, we obtain two disjoint feasible cycles $C^{\prime}$ and $C^{\prime\prime}$ with  $$l(C^{\prime})+l(C^{\prime\prime})\leq l(L_1)+l(L_2)+2<l(C_1)+l(C_2),$$  contradicting the minimality of $l(C_1)+l(C_2)$.\qed

If $f=x_iy_i$ is an edge of $M(C_1)\cup M_0(C_2)$, we use the similar proof of (\ref{eq1}) to consider $C_a-x_iy_i+x_i^-y_i^+$ and $C_b$, where $\{C_a,C_b\}=\{C_1,C_2\}$. Then the followings hold:
\begin{eqnarray}
&\mbox{ If } &t-1\geq s \geq 3, \mbox{ then } e(f,C_1)>\frac{3}{4}l(C_1) \mbox{ for all } f\in M_0(C_2).\label{eq2}\\
&\mbox{ If } &t-1> s =2, \mbox{ then } e(f,C_1)>\frac{3}{4}l(C_1) \mbox{ for all } f\in M_0(C_2).\label{eq3}\\
&\mbox{ If } &s\geq 3, \mbox{ then } e(f,C_2)>\frac{3}{2}t \mbox{ for all } f\in M_0(C_1).\label{eq4}\\
&\mbox{ }  &e(f,C_2)>\frac{3}{2}t \mbox{ for all } f\in M_1(C_1).\label{eq5}
\end{eqnarray}
Next, we claim that
\begin{eqnarray}\label{eq6}s=2 \mbox{ and } t=3.\end{eqnarray}
\textit{\textbf{Proof of (\ref{eq6}).}} If $s\geq 3$, then by (\ref{eq4}), for each $f\in M_0(C_1)$, there exist $f_1,f_2 \in M_0(C_2)$ such that $e(f,f_1)=e(f,f_2)=2$. Otherwise $e(f,C_2)\leq t+1$, contradicting $t\geq 3$.  Thus $[V(C_1\cup C_2)]$ contains two disjoint feasible cycles of order $4$, a contradiction. Thus $s=2$.  Due to Fact \ref{fact1} and $s=2$, we see that $l_{M_1}(C_1)\leq 2$, that is, $l(C_1)=4,6$ or $8$.

If $t\geq 4$, then $t-1\geq 3$. First consider the case of $l(C_1)=4$. Then $e(f,C_1)\geq 4$ for each $f\in M_0(C_2)$ due to (\ref{eq3}) and it follows that $[V(C_1\cup C_2)]$ contains two disjoint feasible cycles of order $4$, a contradiction again.

Therefore, $M_1(C_1)\neq\emptyset$ (i.e. $l_{M_1}(C_1)\geq 1$) and say $f^{\prime}\in M_1(C_1)$. From (\ref{eq5}), $e(f^{\prime},\{y_i,x_{i+2}\})=2$ for some $x_iy_i\in M_0(C_2)$, otherwise $e(f^{\prime},C_2)=\sum_{i=1}^te(f^{\prime},\{y_i,x_{i+2}\})\leq t$. Thus $C_2+f^{\prime}-x_{i+1}y_{i+1}$ has a feasible cycle $C^{\prime}_2$. Then $[\{x_{i+1},y_{i+1}\}\cup V(f)]$ can not be a feasible cycle for all $f\in M_0(C_1)$. So $e(x_{i+1}y_{i+1},C_1)\leq 2+2l_{M_1}(C_1)\leq\frac{3}{4}(4+2l_{M_1}(C_1))=\frac{3}{4}l(C_1)$, which contradicts (\ref{eq2}). Thus $t=3$.\qed

As $\sum_{f\in M_0(C_2)}e(f,C_1)>\frac{3}{4}tl(C_1)$, there exists $f_0\in M_0(C_2)$ such that $e(f_0,C_1)>\frac{3}{4}l(C_1)$. If $l(C_1)=4$, then $e(f_0,C_1)\geq4$.  Furthermore, $e(f_1,C_2)\geq \frac{3}{2}t$ for some $f_1\in M_0(C_1)$ as $e(C_1,C_2)>3t$. Thus $e(f_1,f_2)\geq 2$ for some $f_2\in M_0(C_2)\backslash \{f_0\}$. So $[V(C_1)\cup\{f_0,f_2\}]$ contains two disjoint feasible cycles of order $4$, a contradiction.

Therefore $l(C_1)=6$ or $8$. Let $C_2=x_1y_1\cdots x_3y_3x_1$, where $f_0=x_1y_1$. By (\ref{eq5}), $e(f,C_2)\geq 5$ for all $f\in M_1(C_1)$. Since $e(f_0,C_1)>\frac{3}{4}l(C_1)=\frac{3}{4}(4+2l_{M_1}(C_1))$, there exists $f_1\in M_0(C_1)$ such that $e(f_0,f_1)\geq 2$. Thus $[V(f_0\cup f_1)]$ contains a feasible cycle. Then $C_2-f_0+f$  has no feasible cycle for each $f\in M_1(C_1)$ and so $e(f,\{x_2,y_3\})\leq 1$. Therefore, $$e(f,\{x_1,y_1,x_3,y_2\})=4\mbox{ and } e(f,\{x_2,y_3\})=1,  \mbox{ for each } f\in M_1(C_1).$$ Note that $[V(f\cup f_0)\cup \{x_i,y_i\}]$ ($i=2$ or $i=3$) is a feasible cycle. Then for each $xy\in M_0(C_1)$ and $x_iy_i\in M_0(C_2)$ with $2\leq i\leq 3$, we have $e(\{x,y\},\{x_i,y_i\})\leq 1$, that is, $e(M_0(C_1),\{x_2,y_2,x_3,y_3\})\leq 4$. Hence, $e(C_1,\{x_2,y_2,x_3,y_3\})\leq 4+3l_{M_1}(C_1)$ and it follows immediately that
$$e(C_1,\{x_1,y_1\})>\frac{3}{4}\cdot3\cdot(4+2l_{M_1}(C_1))-(4+3l_{M_1}(C_1))=5+\frac{3}{2}l_{M_1}(C_1).$$
On the other hand, $e(C_1,\{x_1,y_1\})\leq l(C_1)=4+2l_{M_1}(C_1)$, contradicting $l_{M_1}(C_1)\leq 2$.\qed

The following Lemma \ref{lem4} is a key lemma, which is also a directed version of the result of Ore \cite{Ore1960}: Let $u,v$ be the end-vertices of a Hamiltonian path in a graph $G$. If $d_G(u)+d_G(v)\geq |G|$, then $G$ is Hamiltonian. To prove Lemma \ref{lem4}, we need the following Lemma \ref{lem3}. Note that a bipartite tournament is an orientation of a complete bipartite graph and an oriented graph is a digraph without a cycle of order $2$.

\begin{lem}\label{lem3}
If $B=(X,Y)$ is a bipartite oriented graph with $|X|=r_1$ and $|Y|=r_2$, then there exists an arc $uv\in A(B)$ such that $d^-_B(u)+d^+_B(v)<\frac{r}{2}$, where $r=r_1+r_2$.\end{lem}

\proof Suppose to the contrary, let $B$ be the bipartite oriented graph with order $r$ such that for each arc $uv\in A(B)$ we have $d^-_B(u)+d^+_B(v)\geq \frac{r}{2}$. Take the sum of this inequality over all arcs,  we get $$\sum_{uv\in A(B)}(d^-_B(u)+d^+_B(v))\geq\frac{r}{2}|A(B)|.$$

Let $|P_2|$ be the number of paths of length 2 in $B$. Using double counting, one can see that the left side of the above inequality is equal to $2|P_2|$ which implies that
\begin{eqnarray}\label{eq111}
|P_2|\geq\frac{r}{4}|A(B)|.
\end{eqnarray}
Now suppose that $B$ is maximal subject to the above inequality, which is to say, there does not exist an oriented bipartite graph $B^{\prime}$ defined on $(X,Y)$ such that $B^{\prime}$ satisfies the above inequality and $A(B)$ is a proper subset of $A(B^{\prime})$. We now characterize the structure of $B$:

\begin{clm}\label{lemclm1}
For every vertex $u\in V(B)$ with $d^+_B(u)\geq \frac{r}{4}$ or $d^-_B(u)\geq \frac{r}{4}$, $u$ is adjacent to every vertex in the other part.\end{clm}

\proof Without loss of generality, suppose $u\in X$ and there exists a  vertex $v\in Y$ such that $uv\notin A(B)$. We add a new arc $e_0$ to $B$ as follows to get a new graph $B^{\prime}$. If $d^+_B(u)\geq \frac{r}{4}$ we let $e_0$ be $vu$ and,  if $d^-_B(u)\geq \frac{r}{4}$ we let $e_0$ be $uv$. Let $|P_2^{\prime}|$ be the number of  paths with length 2  in $B^{\prime}$ and hence we get $|P_2^{\prime}|\geq |P_2|+\frac{r}{4}$. Thus we have
$$|P_2^{\prime}|\geq |P_2|+\frac{r}{4}\geq\frac{r}{4}|A(B)|+\frac{r}{4}=\frac{r}{4}|A(B^{\prime})|,$$
which contradicts the maximality of $B$.\qed

\begin{clm}\label{lemclm2} There exists a partition of $V(B)$ into $X_1,X_2,Y_1,Y_2$ such that $X=X_1\cup X_2$ and $Y=Y_1\cup Y_2$ such that $(X_1,Y_1)$, $(X_1,Y_2)$ and $(X_2,Y_1)$ are bipartite tournaments and $e(X_2,Y_2)=\emptyset$.
\end{clm}

\proof Recall that we assume that $d^-_B(u)+d^+_B(v)> \frac{r}{2}$ for each arc $uv\in A(B)$.
Then Claim \ref{lemclm1} implies Claim \ref{lemclm2}  by letting $X_1=\{x\in X: d^+_B(x)\geq\frac{r}{4} \mbox{ or } d^-_B(x)\geq\frac{r}{4} \}$, $Y_1=\{y\in Y: d^+_B(y)\geq\frac{r}{4} \mbox{ or } d^-_B(y)\geq\frac{r}{4} \}$, $X_2=X\backslash X_1$ and $Y_2=Y\backslash Y_1$.\qed

Now suppose  $|X_1|=a_1, |X_2|=a_2, |Y_1|=b_1$ and $|Y_2|=b_2$. Clearly, $a_1+a_2=r_1, b_1+b_2=r_2$ and $r_1+r_2=r$. By Claim \ref{lemclm2}, we get that
$$\begin{aligned}
  |P_2|&=\sum_{v\in V(B)}d^+_B(v)d^-_B(v)\\
  &\leq \sum_{v\in V(B)}\frac{(d^+_B(v)+d^-_B(v))^2}{4}\\
&=\sum_{v\in X_1}\frac{(d^+_B(v)+d^-_B(v))^2}{4}+\sum_{v\in Y_1}\frac{(d^+_B(v)+d^-_B(v))^2}{4}\\
& {\,\,\,\,\,}+\sum_{v\in X_2}\frac{(d^+_B(v)+d^-_B(v))^2}{4}+\sum_{v\in Y_2}\frac{(d^+_B(v)+d^-_B(v))^2}{4}\\
&\leq a_1\frac{r_2^2}{4}+b_1\frac{r_1^2}{4}+a_2\frac{b_1^2}{4}+b_2\frac{a_1^2}{4}\\
&\leq \frac{r_1+r_2}{4}(a_1b_1+a_1b_2+a_2b_1)= \frac{r}{4}|A(B)|.\end{aligned}$$

The last inequality can be seen from $a_2b_1^2+b_2a_1^2 \leq a_2b_1r_2+b_2a_1r_1$, $b_1r_1^2 = b_1r_1 a_1+b_1r_1 a_2$ and $a_1r_2^2= a_1r_2b_1+a_1r_2b_2$. This contradicts (\ref{eq111}) and then the lemma follows.\qed

\begin{lem}\label{lem4} Let $P=u_1v_1\cdots u_{r}v_r$ be a good feasible path with $u_i\in X,v_i\in Y$ and $u_iv_{i}\in M(P)$ for all $1\leq i\leq r$. Suppose that $e(\{u_1,v_r\},P)\geq\frac{3}{2}r$ and $[V(P)]$ contains no feasible path $P^{\prime}$ with $M_0(P^{\prime})=M_0(P)$ such that $l(P^{\prime})<l(P)$. Then $[V(P)]$ contains a feasible cycle $C$ such that $V(P)=V(C)$, and then $M_0(C)=M_0(P)$.\end{lem}

\begin{figure}[H]
\centering
\subfigure[$C_{1}\cup C_{2}$] 
{
    \begin{minipage}{5cm}
    \centering          
    \includegraphics[scale=0.4]{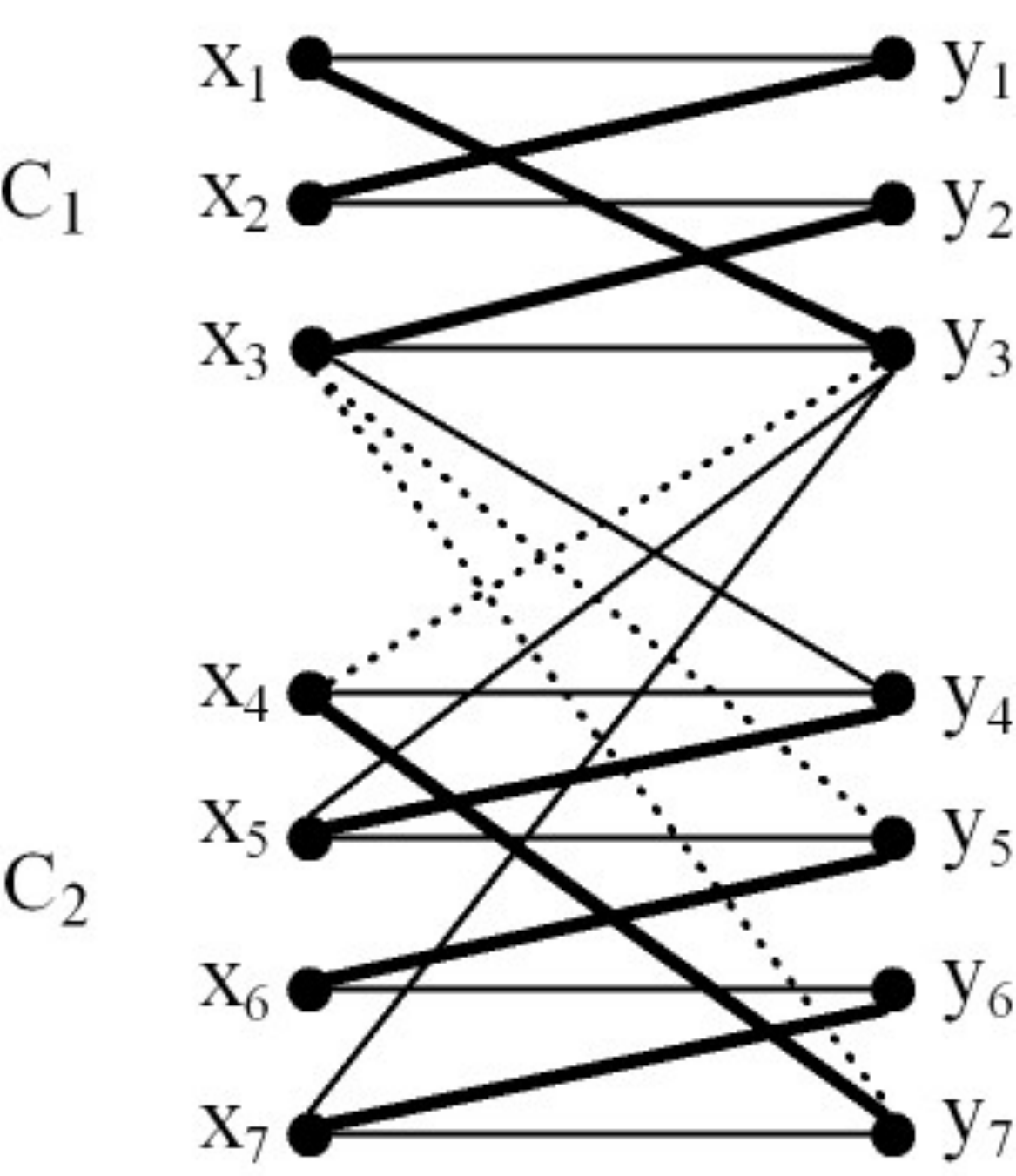}   
    \end{minipage}
}
\subfigure[$B(X,Y)$] 
{
    \begin{minipage}{5cm}
    \centering      
    \includegraphics[scale=0.4]{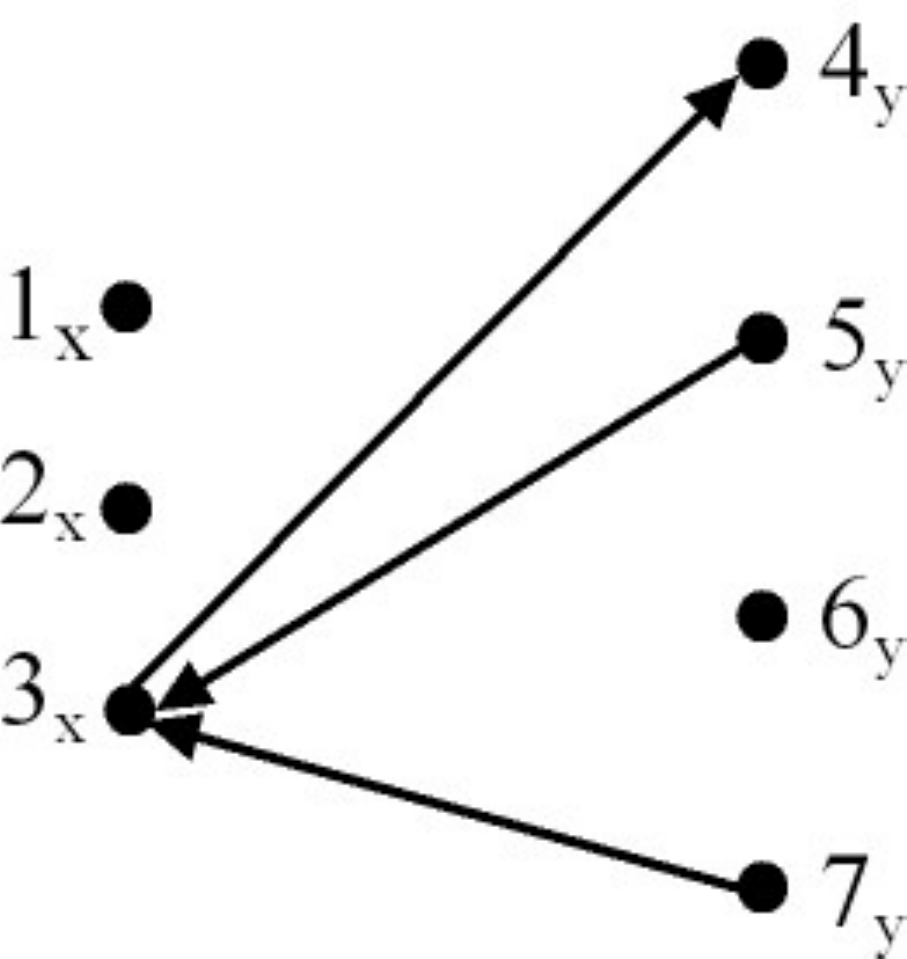}   
    \end{minipage}
}
\caption{An example of $C_{1}\cup C_{2}$ and $B(X,Y)$.} 
\label{fig2}  
\end{figure}

\proof On the contrary, suppose that $[V(P)]$ contains no feasible cycle $C$ such that $V(P)=V(C)$. Since $e(\{u_1,v_r\},P)\geq\frac{3}{2}r$, there exist two vertices $v_i$ and $u_{i+1}$ ($2\leq i\leq r-2$) such that $u_1v_i,u_{i+1}v_r\in E(G)$. Let $C_{1}=u_1v_1\cdots u_{r_1}v_{r_1}u_1$ and $C_{2}=u_{r_1+1}v_{r_1+1}\cdots u_rv_ru_{r_1+1}$, where $r_1+r_2=r$. For convenience, in the following, we relabel $C_{1}=x_1y_1\cdots x_{r_1}y_{r_1}x_1$ and $C_{2}=x_{r_1+1}y_{r_1+1}\cdots x_ry_rx_{r_1+1}$, where $x_i\in X, y_i\in Y$ and $x_iy_i\in E(P)\backslash M(P)$ for each $i$. (See Fig. \ref{fig2} (a).)

Note that since $[V(P)]$ has no feasible cycle $C$ such that $V(P)=V(C)$, we have the following observation: Let $x_p\in V(C_{1})$ and $y_q\in V(C_{2})$, then
\begin{equation}\label{eq15}
\mbox{if }x_py_q\in E, \mbox{ then }x_{q}y_{p}\notin E.
\end{equation}

Let $H$ be the graph induced by the arcs between $C_{1}$ and $C_{2}$ in $[V(P)]$. We define the graph $E_0$ as follows. $V(E_0)=V(H)$ and $E(E_0)=\{x_{q}y_{p}\,|\,x_py_q\in E(H)\}$.

Now construct a bipartite digraph $B(X,Y)$ from $C_{1}$ and $C_{2}$ in the following way.  Let $X=\{1_x,\ldots, (r_1)_x\}$ and $Y=\{(r_1+1)_y,\ldots, r_y\}$. The arc set $A(B)$ of  $B(X,Y)$ is defined by
$$\begin{aligned}
  A(B)&=\{p_xq_y\,|\, x_p\in V(C_1),y_q\in V(C_{2}),x_py_q\in E(H)\}\\
  &\cup\{q_yp_x\,|\,  x_p\in V(C_{1}),y_q\in V(C_{2}), x_py_q\in E(E_0)\}.
\end{aligned}$$

$\bullet$ Fig. \ref{fig2} shows an example with $r_1=3,r_2=4$, where the heavy lines are in $M(P)$, $E(H)=\{x_3y_4,y_3x_5,y_3x_7\}$ and  $E(E_0)=\{y_3x_4,x_3y_5,x_3y_7\}$.
Now we continue to prove Lemma \ref{lem4}. Clearly, $B(X,Y)$ is an oriented graph due to (\ref{eq15}). Furthermore, for $p_x,q_y\in V(B)$, the following hold, where $x_p\in V(C_1), y_q\in V(C_2)$.
$$\begin{aligned}
  &d^+_B(p_x)=d_H(x_p)=d_{E_0}(y_p),d^-_B(p_x)=d_{E_0}(x_p)=d_H(y_p),\\
  &d^+_B(q_y)=d_{E_0}(y_q)=d_H(x_q), \mbox{ and } d^-_B(q_y)=d_{H}(y_q)=d_{E_0}(x_q).
\end{aligned}$$

For arbitrary $x_py_q\in E(E_0)$, one can check that $x_p,y_q$ are the end-vertices of a feasible path $P^{\prime}$ in $[V(P)]$ with $V(P^{\prime})=V(P)$ as $[V(P)]$ contains no shorter feasible path of $M_0$-length $l_{M_0}(P)$. By the definition of good feasible path, we have that $e(\{x_p,y_q\},P)\geq 3r/2$ for each $x_py_q\in E(E_0)$.

Thus for any $p_xq_y\in A(B)$, we have $x_py_q\in E(H)$ ($x_p\in V(C_1), y_q\in V(C_2)$) and then $x_qy_p\in E(E_0)$ due to (\ref{eq15}). Therefore, $$\begin{aligned}
  d^-_B(p_x)+d^+_B(q_y)&=e(\{x_q,y_p\},H)\\
  &=e(\{x_q,y_p\},P)-e(x_q,C_1)-e(y_q,C_2)\\
  &\geq3r/2-r_1-r_2=r/2.
\end{aligned}$$ By symmetry, for each $q_yp_x\in A(B)$, we have $d^-_B(q_y)+d^+_B(p_x)=e(\{x_p,y_q\},H)\geq r/2$, where $x_p\in V(C_1),y_q\in V(C_2)$. This contradicts Lemma \ref{lem3}.\qed

\begin{lem}\label{lem5} Let $H$ be a subgraph of $G$ and $P$ a feasible path with end-vertices $u,v$. Suppose that $P\in G-H$, or $P=uv$ and $P\in H$. Let $f=u^{\prime}v^{\prime}$ be an arbitrary edge of $M_0(G)$ and set $R=\{u,v,u^{\prime},v^{\prime}\}$. If there is no feasible cycle $C$ of $M_0$-length $l_{M_0}(P)+1$ with $V(P\cup f)\subseteq V(C)$, then the following hold.

(i) $e(R,M_1(H))\leq \mbox{max}\{4,3|M_1(H)|\}$. In particular, if $e(R,M_1(H))=4$, then $|M_1(H)|=1$, $e(\{u,v\},f)=0$ and no $f^{\prime}\in M_1(G-(H\cup P))$ linking one vertex of $\{u,v\}$ and one vertex of $\{u^{\prime},v^{\prime}\}$.

(ii) Suppose $f\in H$. Then $e(R,M_0(H))\leq e(\{u,v\},f)+3|M_0(H)|-1$. In particular, if $P=uv\in H$ ($l_{M_0(P)}=1$), then$$e(R,M_0(H))\leq\left\{\begin{aligned}
  &2|M_0(H)|,  \,\,\,e(f,P)=0\\
  &2|M_0(H)|+2,  \,\,\,e(f,P)=1.
\end{aligned}
\right.$$
\end{lem}
\proof(i) If $|M_1(H)|=1$, then (i) holds clearly. Thus assume that $|M_1(H)|\geq 2$. Suppose (i) is not true, then  there are two edges $f_1$ and $f_2$ in $M_1(H)$ such that $e(R,\{f_1,f_2\})\geq 7$. One can check easily that $[V(P\cup f_1\cup f_2\cup f)]$ contains a feasible cycle $C$ of $M_0$-length $l_{M_0}(P)+1$ with $V(P\cup f)\subseteq V(C)$, a contradiction.

(ii) Since $H$ contains no feasible cycle of length $M_0$-length $l_{M_0}(P)+1$, then $e(\{u,v\},f^{\prime})\leq 1$ for each $f^{\prime}\in M_0(H)$. Clearly, $e(f_1,f_2)\leq 2$ for all $f_1,f_2\in M_0(H)$. Thus $e(R,M_0(H))\leq e(\{u,v\},f)+|M_0(H)|-1+2|M_0|=e(\{u,v\},f)+3|M_0(H)|-1$. Furthermore, if $P=uv\in H$, we have that $e(f_1,f_2)\leq 1$ for all $f_1,f_2\in M_0(H)$. Thus $e(f,P)\leq1$ and $e(R,M_0(H))\leq e(R)+2(|M_0(H)|-2)=e(R)+2|M_0(H)|-4$, and so (ii) holds.\qed

\vspace{-6mm}
\section{Proof of Theorem \ref{thm2}}  \label{section 5}
Let $G(X,Y)$ be a balanced bipartite graph with perfect matching $M=\{x_iy_i: x_i\in X, y_i\in Y, 1\leq i\leq n\}$. Suppose that $M_0\subseteq M$ and $\delta(M_0)\geq (3n+1)/4$. To prove Theorem \ref{thm2}, we prove the following claim first.

\begin{clm}\label{clm3} There exist $\lfloor\frac{|M_0|}{2}\rfloor$ disjoint feasible cycles $C_1,\ldots, C_{\lfloor\frac{|M_0|}{2}\rfloor}$ with $l_{M_0}(C_i)=2$ for each $i$.
\end{clm}
\proof Suppose, for a contradiction, that the claim  fails. We may assume that $n\geq2$ and $|M_0|\geq 2$. Let $p$ be the largest number such that $G$ contains $p$ disjoint feasible cycles $C_1,\ldots, C_p$. We claim that $p\geq 1$.  Let $f_1=uv,f_2=u^{\prime}v^{\prime}\in M_0$ and $R=\{u,v,u^{\prime},v^{\prime}\}$. If $p=0$, then from Lemma \ref{lem5} (let $H=G$), we see that $e(R,G)=e(R,M_0)+e(R,M_1)\leq \mbox{max}\{2|M_0|+4,2|M_0|+2+3|M_1|\}$. However, $e(R,G)\geq 4\cdot\frac{3n+1}{4}$, where $n=|M_0|+|M_1|$, a contradiction. So $p\geq 1$.
\medskip

We choose $p$ disjoint feasible cycles $C_1,\ldots, C_p$ such that

(\MyRoman{1}) $\sum_{i=1}^{p}l_{M_0}(C_i)$ is minimum, and

(\MyRoman{2}) $\sum_{i=1}^{p}l(C_i)$ is minimum, subject to (\MyRoman{1}).

\medskip
Say $l_{M_0}(C_1)\leq \cdots \leq l_{M_0}(C_p)$. Set $\mathcal{C}=\cup_{i=1}^{p}C_i$ and  $H=G-\mathcal{C}$. Furthermore, $2c,2h,m_0$ and $m_1$, respectively, denote  $|V(\mathcal{C})|, |V(H)|, |M_0(H)|$ and $|M_1(H)|$. So, $c=n-h$ and $h=m_0+m_1$.

Next we claim that $p=\lfloor\frac{|M_0|}{2}\rfloor$ and $l_{M_0}(C_p)=2$. To see this, we first suppose, for a contradiction, that $l_{M_0}(C_p)=t\geq 3$. Let $C_p=x_1y_1\cdots x_{l}y_{l}x_1$, where $l=l(C_p)/2$ and $f_i=x_iy_i\in M(C_p)$. By the choice (\MyRoman{1}), for each $f_i\in M_0(C_p)$ and $f\in M(C_p)\cup M_0(H)$, we conclude that
$$e(f_i,f)\leq\left\{\begin{aligned}
  &2,  \,\,\,f\in M_1(C_p)\cap\{f_{i+1},f_{i-1}\}\\
  &1,  \,\,\,f\in M_0(H)\cup(M_0(C_p)\cap\{f_{i+1},f_{i-1}\})\\
  &0,  \,\,\,\mbox{Others}.
\end{aligned}
\right.$$

Thus $\sum_{f_i\in M_0(C_p)}e(f_i,M_0(H))\leq m_0t.$ In addition, we claim that $\sum_{f_i\in M_0(C_p)}e(f_i,C_p)\leq\frac{3}{4}tl(C_p)$ and $\sum_{f_i\in M_0(C_p)}e(f_i,M_1(H))\leq \frac{3}{2}tm_1$. For each $f_i\in M_0(C_p)$, let $\lambda_i=|M_1(C_p)\cap\{f_{i-1},f_{i+1}\}|$, then $e(f_i,C_p)\leq 2+2\lambda_i+(2-\lambda_i)=4+\lambda_i$. Note that $\sum_{f_i\in M_0(C_p)}\lambda_i=2|M_1(C_p)|$, $l(C_p)=2t+2|M_1(C_p)|$ and $t\geq 3$. It follows that $$\sum_{f_i\in M_0(C_p)}e(f_i,C_p)\leq 4t+\sum_{f_i\in M_0(C_p)}\lambda_i=l(C_p)+2t\leq 2l(C_p)<\frac{3}{4}tl(C_p).$$

Moreover, let $f_i=x_iy_i, f_j=x_jy_j$ be two independent edges of $M_0(C_p)$ with $f_j$ being  next to $f_i$ in $M_0(C_p)$. Note that $e(f_i,f_j)=1$ or there is an edge of $M_1(C_p)$ linking $f_i$ and $f_j$. Combining Lemma \ref{lem5} (i) and the choice (\MyRoman{1}), we have  that $e(\{x_i,x_j,y_i,y_j\},M_1(H))\leq 3m_1$ and then $\sum_{f_i\in M_0(C_p)}e(f_i,M_1(H))\leq \frac{3}{2}tm_1.$ Therefore,\begin{eqnarray}
\sum_{f_i\in M_0(C_p)}e(f_i,\mathcal{C}-C_p)&=&2t\cdot\frac{3n+1}{4}-\sum_{f_i\in M_0(C_p)}(e(f_i,M_0(H))+e(f_i,C_p)+e(f_i,M_1(H)))\nonumber\\
&>&2t\cdot\frac{3n+1}{4}-m_0t-\frac{3}{4}tl(C_p)-\frac{3}{2}tm_1\nonumber\\
&=&\frac{3t}{4}(2c-l(C_p))+\frac{t+m_0t}{2}\nonumber\\
&>&\frac{3t}{4}(2c-l(C_p))\nonumber.\end{eqnarray} This implies that $\sum_{f_i\in M_0(C_p)}e(f_i,C_j)>\frac{3}{4}tl(C_j)$ for some $1\leq j\leq p-1$. By Lemma \ref{lem2}, we obtain a contradiction to  the choice (\MyRoman{1}) or (\MyRoman{2}). Thus $l_{M_0}(C_i)=2$ for all $1\leq i\leq p$. According to Fact \ref{fact1}, we have that $l(C_i)\leq 8$ for all $1\leq i\leq p$.

Now we show that $p=\lfloor\frac{|M_0|}{2}\rfloor$. Assume that $p<\lfloor\frac{|M_0|}{2}\rfloor$, then $m_0\geq 2$. Let $f_1=x_1y_1$ and $f_2=x_2y_2$ be two distinct edges of $M_0(H)$ and let $R=\{x_1,y_1,x_2,y_2\}$. Due to Lemma \ref{lem5}, we conclude that $e(R,H)\leq \mbox{max}\{2m_0+4,2m_0+2+3m_1\}=2m_0+2+3m_1$. Hence $e(R,\mathcal{C})\geq4\cdot(3n+1)/4-(2m_0+3m_1+2)\geq3c+1$. This means that there is a $C_i$ with $1\leq i\leq p$ such that $e(R,C_i)\geq \frac{3}{2}l(C_i)+1$.

Recall that $l_{M_0}(C_i)=2$ and $l(C_i)\leq 8$. Say $M_0(C_i)=\{xy,uv\}$. If $l(C_i)=4$, then $e(R,C_i)\geq 7$. Then each of $[\{x_1,y_1,x,y\}]$ and $[\{x_2,y_2,u,v\}]$ contains a feasible cycle, contradicting that $p$ is maximum. Thus $l(C_i)=6$ or $l(C_i)=8$.  Note that $1\leq l_{M_1}(C_i)\leq 2$. Choices (\MyRoman{1}) and (\MyRoman{2}) imply that $e(R,f)\leq 4$ for each $f\in M_1(C_i)$ and $e(R,f)\leq 2$ for each $f\in M_0(C_i)$. So $e(R,C_i)\leq 4+4l_{M_1}(C_i)$. However, $e(R,C_i)\geq\frac{3}{2}l(C_i)+1=\frac{3}{2}(4+2l_{M_1}(C_i))+1$, which contradicts $l_{M_1}(C_i)\leq 2$.\qed

\noindent\textbf{Proof of Theorem \ref{thm2}. } Now we are ready to prove Theorem \ref{thm2}.  By way of contradiction, suppose Theorem \ref{thm2} fails. That is, $G$ does not contain $k$ disjoint cycles of $M_0$-length $n_1,\ldots,n_k$, respectively. Among all such integers, we choose $n_i$ ($1\leq i\leq k$) with $n_1+\cdots+n_k$ minimal.

The case $n_1=\cdots=n_k=2$ is proved by Claim \ref{clm3}. Thus we may assume that $n_1+\cdots+n_k>2k$. By the minimality of $\sum_{i=1}^kn_i$, for some $1\leq i\leq k$  with $n_i\geq 3$, $G$ contains $k$ disjoint feasible cycles of $M_0$-lengths $n_1,\ldots,n_{i-1},n_i-1,n_{i+1},\ldots, n_k$, respectively. Therefore, for some $1\leq j\leq k$, $G$ contains $k-1$ disjoint feasible cycles $C_1,\ldots, C_{j-1},C_{j+1},\ldots, C_k$ of $M_0$-lengths $n_1,\ldots,n_{j-1},n_{j+1},\ldots, n_k$, respectively and
a good feasible path $P$ of $M_0$-length $n_j-1$ such that $P$ is disjoint from all these cycles. These $k-1$ cycles and the path $P$ exist, since for any $M$-alternating path $P^{\prime}$, there is a good feasible path $P$ in $[V(P^{\prime})]$. We choose $C_1,\ldots, C_{j-1},C_{j+1},\ldots, C_k$ and $P$ such that
\begin{equation}\label{eq9}
l(P)+\sum_{i\neq j}l(C_i) \mbox{ is minimum}\tag{$\ast$}.
\end{equation}

For convenience, say $l_{M_0}(C_i)=n_i$ for each $1\leq i\leq k-1$, $l_{M_0}(P)=n_k-1$, and $P=x_1y_1\cdots x_ry_r$. Let $\mathcal{C}=\cup_{i=1}^{k-1}C_i$ and $H=G-(\mathcal{C}\cup P)$. Furthermore, $2c,2h,m_0$ and $m_1$, respectively, denote  $|V(\mathcal{C})|, |V(H)|, |M_0(H)|$ and $|M_1(H)|$. So, $r=n-h-c$ and $h=m_0+m_1\geq 1$.  Let $xy\in M_0(H)$ and $R=\{x_1,y_r,x,y\}$.
\begin{clm}\label{clm4} $e(R,H)\leq 3h$ and $e(R,C_j)\leq\frac{3}{2}l(C_j)$ for all $1\leq j\leq k-1$.\end{clm}

\proof Since $[V(P\cup H)]$ contains  no feasible cycle of $M_0$-length $n_k$, $e(R,H)\leq 3m_0+3m_1=3h$ due to Lemma \ref{lem5}. Next assume that there exists $C_j$ with $1\leq j\leq k-1$ such that $e(R,C_j)\geq \frac{3}{2}l(C_j)+1$. Let $f_i=x_iy_i$ for $1\leq i\leq n_j$ be a list of edges of $M_0(C_j)$ along the direction of $C_j$. Let $f_i^+$ be the next edge of $f_i$ along the direction of $C_j$ and  $B=\{i\,|\,f_i^+\in M_1(C_j),1\leq i\leq n_j\}$ and $A=\{1,\ldots, n_j\}\backslash B$. Let $$s_i=
\begin{cases}
e(\{y_i,x_{i+2}\},\{x,y\})+e(\{x_{i+1},y_{i+1}\},\{x_1,y_r\}),& i\in A,\\
e(\{y_i,x_{i+2}\},\{x,y\})+e(\{x_{i+1},y_{i+1}\},\{x_1,y_r\})+e(f_i^+,R),& i\in B.
\end{cases}$$
Then \begin{eqnarray}
\frac{3}{2}l(C_j)+1\leq e(R,C_j)=\sum_{i\in A}s_i+\sum_{i\in B}s_i.
\label{eq11}
\end{eqnarray} where the subscripts $i$ are taken modulo in $\{1,\ldots,n_j\}$. On the other hand,
\begin{eqnarray}
\frac{3}{2}l(C_j)&=&\frac{3}{4}\sum_{i\in A}[|\{y_i,x_{i+2}\}|+|\{x_{i+1},y_{i+1}\}|]\nonumber\\
&+&\frac{3}{4}\sum_{i\in B}[|\{y_i,x_{i+2}\}|+|\{x_{i+1},y_{i+1}\}|+2|V(f_i^+)|].\label{eq12}
\end{eqnarray}
Combining (\ref{eq11}) and (\ref{eq12}), we obtain that either there exists $i\in A$ such that $s_i\geq (3/4)\cdot 4+1=4,$ or there exists $i\in B$ such that $s_i\geq (3/4)\cdot (4+4)+1=7.$

If the former holds, we see that $e(\{y_i,x_{i+2}\},\{x,y\})=e(\{x_{i+1},y_{i+1}\},\{x_1,y_r\})=2$. Consequently, $[(V(C_j)\backslash \{x_{i+1},y_{i+1}\})\cup\{x,y\}]$ contains a feasible cycle of $M_0$-length $n_j$ and $[V(P)\cup\{x_{i+1},y_{i+1}\}]$ contains a feasible cycle of $M_0$-length $n_k$, a contradiction. Therefore, the latter holds. As above, we see that $e(\{y_i,x_{i+2}\},\{x,y\})+e(\{x_{i+1},y_{i+1}\,\{x_1,y_r\})=3$ and $e(f_i^+,R)=4$. If $e(\{y_i,x_{i+2}\},\{x,y\})=2$, then $[V(P\cup(C_i\backslash f_i^+))\cup\{x,y\}]$ contains a feasible cycle of $M_0$-length $n_j$ and a path $P$ with $l_{M_0}(P)=n_k-1$. This contradicts the choice (\ref{eq9}). The case of $e(\{x_{i+1},y_{i+1}\},\{x_1,y_r\})=2$ is similar. Thus $e(R,C_j)\leq\frac{3}{2}l(C_j)$ for all $1\leq j\leq k-1$.\qed

Therefore, $e(R,P)\geq 4\cdot\frac{3n+1}{4}-3c-3h=3r+1$ holds by Claim \ref{clm4}. Recall that $|P|=2r$. If $r=1$, i.e. $n_k=2$. Since $e(R,P)\geq 4$, $[V(R)]$ contains a feasible cycle $C$ of $l_{M_0}(C)=2=n_k$, a contradiction. Thus $r\geq 2$.

\bigskip
\noindent\textbf{Case 1. } $e(\{x_1,y_r\},P)\geq3r/2$.

Note that $e(\{x_1,y_r\},P)\leq 2r$  and then $e(\{x,y\},P)\geq r+1$. From the choice (\ref{eq9}), $[V(P)]$ contains no feasible path $P^{\prime}$ with $M_0(P^{\prime})=M_0(P)$ such that $l(P^{\prime})<l(P)$. Moreover, $[V(P)]$  contains  no feasible cycle of $M_0$-length $n_k-1$, otherwise, $[V(P)\cup\{x,y\}]$ has a feasible cycle of $M_0$-length $n_k$ from Lemma \ref{lem1}, a contradiction. It follows by Lemma \ref{lem4} that $e(\{x_1,y_r\},P)<3r/2$, a contradiction again.

\bigskip

Similarly by the argument with Case 1, we have $x_1y_r\not\in E(G)$, i.e. $e(\{x_1,y_r\},P)\leq 2r-2$ and then $e(\{x,y\},P)\geq r+3$. Furthermore, if there exists $f\in M_1(H)$ such that $[V(P\cup f)]$ contains a feasible cycle $C$ with $V(P)\subseteq V(C)$ and $l_{M_0}(C)=n_k-1$, then by Lemma \ref{lem1}, $[V(P\cup f)\cup \{x,y\}]$ has a feasible cycle of $M_0$-length $n_k$ as $e(\{x,y\},C)\geq e(\{x,y\},P)\geq r+3\geq l(C)/2+1$, a contradiction again. Thus $e(\{x_1,y_r\},f)\leq 1$ for each  $f\in M_1(H)$. As $[V(P\cup H)]$ has no feasible cycle of $M_0$-length $n_k$, we know that $e(\{x_1,y_r\},f)\leq 1$ for every $f\in M_0(H)$. This means that $e(\{x_1,y_r\},H)\leq m_0+m_1=h$.

\bigskip

\noindent\textbf{Case 2. } $e(\{x_1,y_r\},P)<3r/2$.

In this case, we see that $e(\{x_1,y_r\},P\cup H)<3r/2+h$ and
\begin{eqnarray}
e(\{x,y\},P)> 3r+1-3r/2=3r/2+1, \mbox{ i.e. }e(\{x,y\},P)\geq 3(r+1)/2.\label{eq13}
\end{eqnarray}

So $e(\{x_1,y_r\},\mathcal{C})> 2\cdot(3n+1)/4-(3r/2+h)=3c/2+(h+1)/2$. This means that there exists $C_j\in \mathcal{C}$ such that \begin{eqnarray} e(\{x_1,y_r\},C_j)> 3l(C_j)/4.\label{eq14}\end{eqnarray}

The definition of $A$ and $B$ is the same as in Claim \ref{clm4}. It follows by Fact \ref{fact1} that $|B|\leq |M_0(C_j)|$. Therefore, \begin{eqnarray}
e(\{x_1,y_r\},M_0(C_j))&>&3l(C_j)/4-e(\{x_1,y_r\},M_1(C_j))\nonumber\\
&\geq&3(2|M_0(C_j)|+2|B|)/4-2|B|\nonumber\\
&\geq&|M_0(C_j)|.\nonumber
\end{eqnarray}
So $e(\{x_1,y_r\},f)\geq2$ for some $f\in M_0(C_j)$ and then $[V(P\cup f)]$ contains a feasible cycle $C^{\prime}$ of $M_0$-length $n_k$ and $C_j\backslash f$ contains a good feasible path $P^{\prime}=x^{\prime}_1y^{\prime}_1\cdots x^{\prime}_{r^{\prime}}y^{\prime}_{r^{\prime}}$ with $l_{M_0}(P^{\prime})=n_j-1$. Clearly, $l(C^{\prime}_j)+l(P^{\prime})\leq l(C_j)+l(P)$. From (\ref{eq9}), we see that $l(C^{\prime}_j)+l(P^{\prime})=l(C_j)+l(P)$. So $l(P^{\prime})=2r^{\prime}=l(C_j)-2$. With $C^{\prime}$ and $P^{\prime}$ in place of $C_j$ and $P$ in the above argument:
if $e(\{x^{\prime}_1,y^{\prime}_{r^{\prime}}\},P^{\prime})\geq3r^{\prime}/2$, then by Case 1, we are done. Thus $e(\{x^{\prime}_1,y^{\prime}_{r^{\prime}}\},P^{\prime})<3r^{\prime}/2$. Furthermore, we obtain $e(\{x,y\},P^{\prime})\geq3(r^{\prime}+1)/2$ as we obtain (\ref{eq13}). It follows by (\ref{eq14}) that
\begin{eqnarray}
e(R,C_j)&\geq& e(\{x_1,y_r\},C_j)+e(\{x,y\},C_j)\nonumber\\
&>&\frac{3}{4}l(C_j)+\frac{3}{2}(r^{\prime}+1)\nonumber\\
&=&\frac{3}{4}l(C_j)+\frac{3}{4}(l(C_j)-2)+\frac{3}{2}\nonumber\\
&=&\frac{3}{2}l(C_j),\nonumber
\end{eqnarray} contradicting Claim \ref{clm4}. This proves the theorem.\qed
\vspace{-6mm}
\section{Proof of Theorems \ref{thmcyclable} and \ref{thm5}}  \label{section 6}
First, we prepare terminology and notation which will be used in the following proofs. Let $C$ be a cycle.  A \emph{generalized C-bypass} $T$ (with respect to $W$) is a path such that only the end-vertices $x$ and $y$ belong to $V(C)$ and at least one internal vertex belongs to $W$.  Note that we allow $x$ and $y$ to be the same vertex. Let $P=v_1v_2\cdots v_p$ be a path and $u$  a vertex not on $P$. If there are two vertices $v_i$ and $v_{i+1}$ such that $v_iu,uv_{i+1}\in A(D)$, then we say that $u$ can be \emph{inserted} into $P$.

\medskip

\noindent\textbf{Proof of Theorem \ref{thmcyclable}. }

Suppose, for a contradiction, that Theorem \ref{thmcyclable} fails. Let $D$ be a digraph of order $n$ and  $W\subseteq V(D)$ with  $\delta^0(W)\geq n/2$. First, we present two claims.

\begin{clm}\label{clm61}
Let $u$ and $v$ be two vertices in $W$. Then there is a path $P_1$ of length at most $2$ from $u$ to $v$ and a path $P_2$ of length at most $2$ from $v$ to $u$. Moreover, $P_1$ and $P_2$ are internally disjoint.\end{clm}
\proof  Since $\delta^0(W)\geq n/2$, we see that $d^+_D(u)+d^-_D(v)\geq n$ and $d^+_D(v)+d^-_D(u)\geq n$. Thus $uv\in A(D)$ or there are at least two vertices in $N^+_D(u)\cap N^-_D(v)$. By symmetry, we have that $vu\in A(D)$ or there are at least two vertices in $N^+_D(v)\cap N^-_D(u)$. Then we can find such $P_1$ and $P_2$ easily.\qed
\vspace{-6mm}
\begin{clm}\label{clm63}(\cite{Berman1998})
Let $P$ be a path of length $p$ in $D$ and $u$ a vertex not on $P$. If $d_P(u)\geq p+2$, then $u$ can be inserted into $P$.\end{clm}

According to Claim \ref{clm61}, there is a cycle $C$ which contains at least two vertices of $W$. Thus we may assume that $|W|\geq 3$. Since $W$ is not cyclable, assume that $v$ is a vertex of $W$ which is not in $C$. From Claim \ref{clm61}, it is easy to verify that there exists a generalized $C$-bypass $T$ with origin $x$ and terminus $y$ such that $v\in V(T)$, where $x,y\in V(C)$. Set $P_1=T[x,v]$, $P_2=T[v,y]$, $S=C[x^+,y^-]$ and $B=C[y,x]$.

We choose a cycle $C$ and a generalized $C$-bypass $T$ such that

(\MyRoman{1}) $C$ contains as many vertices of $W$ as possible,

(\MyRoman{2}) $|S|$ is minimum, and

(\MyRoman{3}) $|T|$ is minimum.
\bigskip

Similar to the notation $M(H)$, we define $W(H)$ to be $W\cap V(H)$. Let $u\in W(C)$ and $R=V(D-C)$.

\begin{clm}\label{clm62}(\cite{Berman1998})
(i) Let $v\in W(T\backslash\{x,y\})$ and $u\in W(S)$. If $d_R(v)+d_R(u)\leq 2|R|-2$ and $d_B(v)\leq |B|+1$, then $u$ can be inserted into $B$.

(ii) Let $P$ and $Q$ be two disjoint paths and $K\subseteq V(P)$. If every vertex $z$ in $K$ can be inserted into $Q$, then there exists a path $Q^{\prime}$ such that $K\cup V(Q)\subseteq V(Q^{\prime})$.\end{clm}

We are now ready to complete the proof of Theorem \ref{thmcyclable}. Let $x_1,\ldots, x_k$ be the sequence of vertices of $W$ on $S$ listed in the order they occur along $C$.

Suppose $w_1\in N^+_{R\backslash P_2}(x_i)\cap N^-_{R\backslash P_2}(v)$ for some $x_i\in W(S)$. Then replacing $T$ with the new generalized $C$-bypass $T^{\prime}=x_iw_1vP_2$ reduces the size of $S$, contradicting the minimality of $S$. Thus we have
\begin{equation}\label{eq61}
N^+_{R\backslash P_2}(x_i)\cap N^-_{R\backslash P_2}(v)=\emptyset \mbox{ for each } x_i\in W(S).
\end{equation}

By symmetry, we also have
\begin{equation}\label{eq62}
N^-_{R\backslash P_1}(x_i)\cap N^+_{R\backslash P_1}(v)=\emptyset \mbox{ for each } x_i\in W(S).
\end{equation}

Clearly, due to the choice (\MyRoman{1}), we may assume that $v$ can not be inserted into $B$. Thus by Claim \ref{clm63}, we have
\begin{equation}\label{eq63}
d_B(v)\leq |B|+1.
\end{equation}

Next, we will show  that\begin{eqnarray}
N^+_{P_2\backslash y}(x_i)\cap N^-_{P_2\backslash y}(v)=\emptyset \mbox{ for each } x_i\in W(S),\label{eq64}\\
N^-_{P_1\backslash x}(x_i)\cap N^+_{P_1\backslash x}(v)=\emptyset \mbox{ for each } x_i\in W(S).\label{eq65}
\end{eqnarray}

By symmetry, it is sufficient to prove (\ref{eq64}).
\medskip

\noindent\textit{\textbf{Proof of (\ref{eq64}). }}  If not, say $N^+_{P_2\backslash y}(x_i)\cap N^-_{P_2\backslash y}(v)\neq\emptyset$  for some $x_i\in W(S)$. Set $t=\mbox{ max}\{i\,|\,N^+_{P_2\backslash y}(x_i)\cap N^-_{P_2\backslash y}(v)\neq\emptyset, 1\leq i\leq k\}$. By the choice of $t$, we have (\ref{eq64}) holds for all $i$ with $t+1\leq i\leq k$. In this case, we have (\ref{eq65}) holds for all $i$ with $t+1\leq i\leq k$. For otherwise, setting $w_1\in N^-_{P_1\backslash x}(x_i)\cap N^+_{P_1\backslash x}(v)$ for some $i\in\{t+1,\ldots, k\}$ and $w_2\in N^+_{P_2\backslash y}(x_t)\cap N^-_{P_2\backslash y}(v)$, we obtain the new generalized $C$-bypass $x_tw_2vw_1x_i$ with $S=C[x^+_t,x^-_i]$, contradicting the choice (\MyRoman{2}).

Therefore, Claims \ref{clm61} and  \ref{clm62} (i) imply that  $N^+_{R}(x_i)\cap N^-_{R}(v)=\emptyset$ and $N^-_{R}(x_i)\cap N^+_{R}(v)=\emptyset$ for each $t+1\leq i\leq k$. And hence, $|N^+_{R}(x_i)|+|N^-_{R}(v)|\leq |R|-1$ and $|N^-_{R}(x_i)|+|N^+_{R}(v)|\leq |R|-1$, that is, $d_R(x_i)+d_R(v)\leq 2|R|-2$ for $i=t+1,\ldots, k$. Applying Claim \ref{clm62} and (\ref{eq63}), there exits a path $B^{\prime}$ such that $\{x_{t+1},\ldots, x_k\}\cup V(B)\subseteq V(B^{\prime})$. Obviously, $C^{\prime}=B^{\prime}C[x,x_t]x_tw_2P_2[w_2,y]$ forms a cycle, where $w_2\in N^+_{P_2\backslash y}(x_t)\cap N^-_{P_2\backslash y}(v)$. Furthermore, $T[x,w_2]$ is a generalized $C^{\prime}$-bypass, contradicting the minimality of $T$ (note that $S$ is not increased). Thus, (\ref{eq64}) holds.\qed

By (\ref{eq61}),(\ref{eq62}),(\ref{eq64}) and (\ref{eq65}), we obtain $N^+_{R}(x_i)\cap N^-_{R}(v)=\emptyset$ and $N^-_{R}(x_i)\cap N^+_{R}(v)=\emptyset$ for each $1\leq i\leq k$, that is, $d_R(x_i)+d_R(v)\leq 2|R|-2$ for all $1\leq i\leq k$. Using Claim \ref{clm62}, there is a path $B^{\prime}$ from $y$ to $x$, where $W(S)\cup V(B)\subseteq V(B^{\prime})$. Then $TB^{\prime}$ forms a cycle containing more vertices of $W$ than $C$ does, contradicting the choice (\MyRoman{1}), and Theorem \ref{thmcyclable} is proved.\qed

\noindent\textbf{Proof of Theorem \ref{thm5}. }
Let $D$ be a digraph of order $n\geq 2|W|$. Suppose that $\delta^0(W)\geq \frac{n}{2}+|W|-1$. We prove that $D$ contains an arbitrary $W$-cycle-factor.

\begin{figure}[H]
    \centering      
    \includegraphics[scale=0.5]{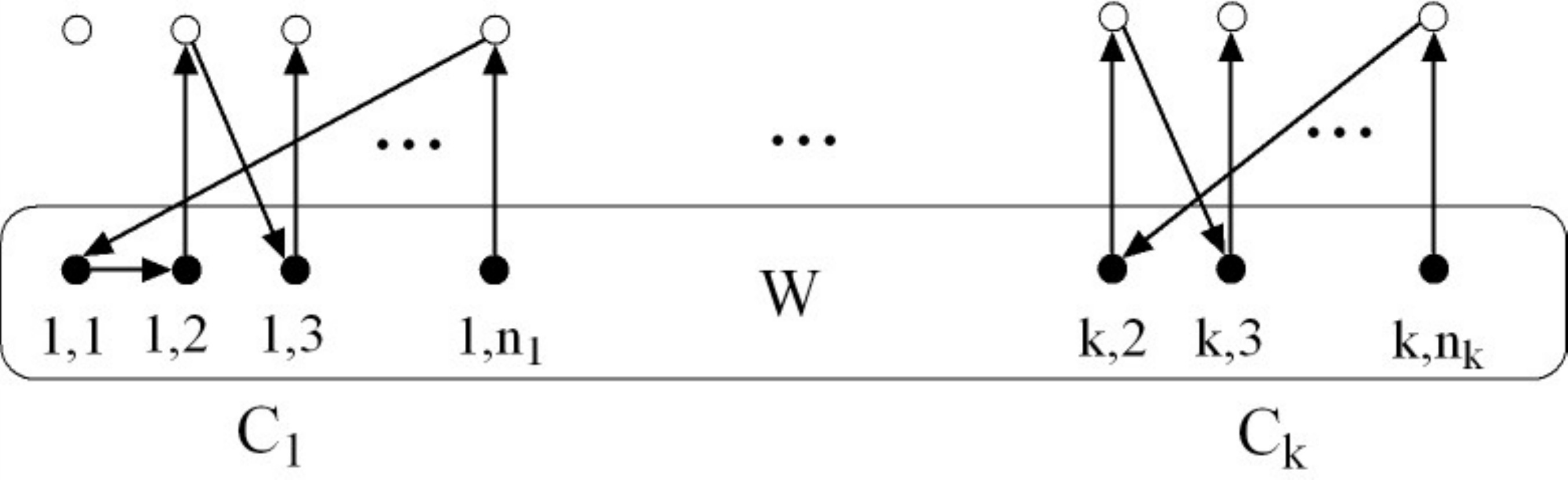}    
\caption{$k$ disjoint cycles with $|V(C_i)\cap W|=n_i$.}
\label{fig3} 
\end{figure}

For any positive integer partition $|W|=n_1+\cdots+n_k$ with $k\geq 1$ and $n_i\geq 2$ for each $i$, let $W=\{w_{1,1},\ldots,w_{1,n_1},\ldots,w_{k,1},\ldots,w_{k,n_k}\}$.  Since $d^+_D(w_{i,j})+d^-_D(w_{i,j+1})\geq n+2|W|-2$ for each $1\leq i \leq k,1\leq j\leq n_i$ (the subscripts $j$ are taken modulo in $\{1,\ldots, n_i\}$), then $w_{i,j}w_{i,j+1}\in A(D)$ or there are at least $|W|$ vertices in $N^+_{D-W}(w_{i,j})\cap N^-_{D-W}(w_{i,j+1})$. In this way, one can easily find an arbitrary $W$-cycle-factor (see Fig. \ref{fig3}, where the vertex $w_{i,j}$ is replaced by $i,j$).\qed

\section{Remark}  \label{section 7}

In this section, we discuss the degree condition in Theorem \ref{thm1}. Note that no cycle in the digraph $D$ contains exactly one vertex of $W$ when $W=V(D)$, so  $n_i\geq 2$ for each $i$ in Theorem \ref{thm1}. Furthermore, the minimum semi-degree condition in Theorem \ref{thm1} is sharp in some sense.

Consider the digraph $D_1$ which consists four parts $U,X,Y$ and $Z$, where $U=X=K^*_{4k-1}$ and $Y=Z=K^*_{4k}$. The bold arcs indicate complete domination in the direction shown. Note that $D_1$ has order $16k-2$ and $\delta^0(D_1)=12k-3<(3n-3)/4$. However, $D_1$ does not contain $8k-1$ disjoint cycles of order 2, as $|U\cup Y|$ is odd.
\begin{figure}[H]
\centering    
    \includegraphics[scale=0.4]{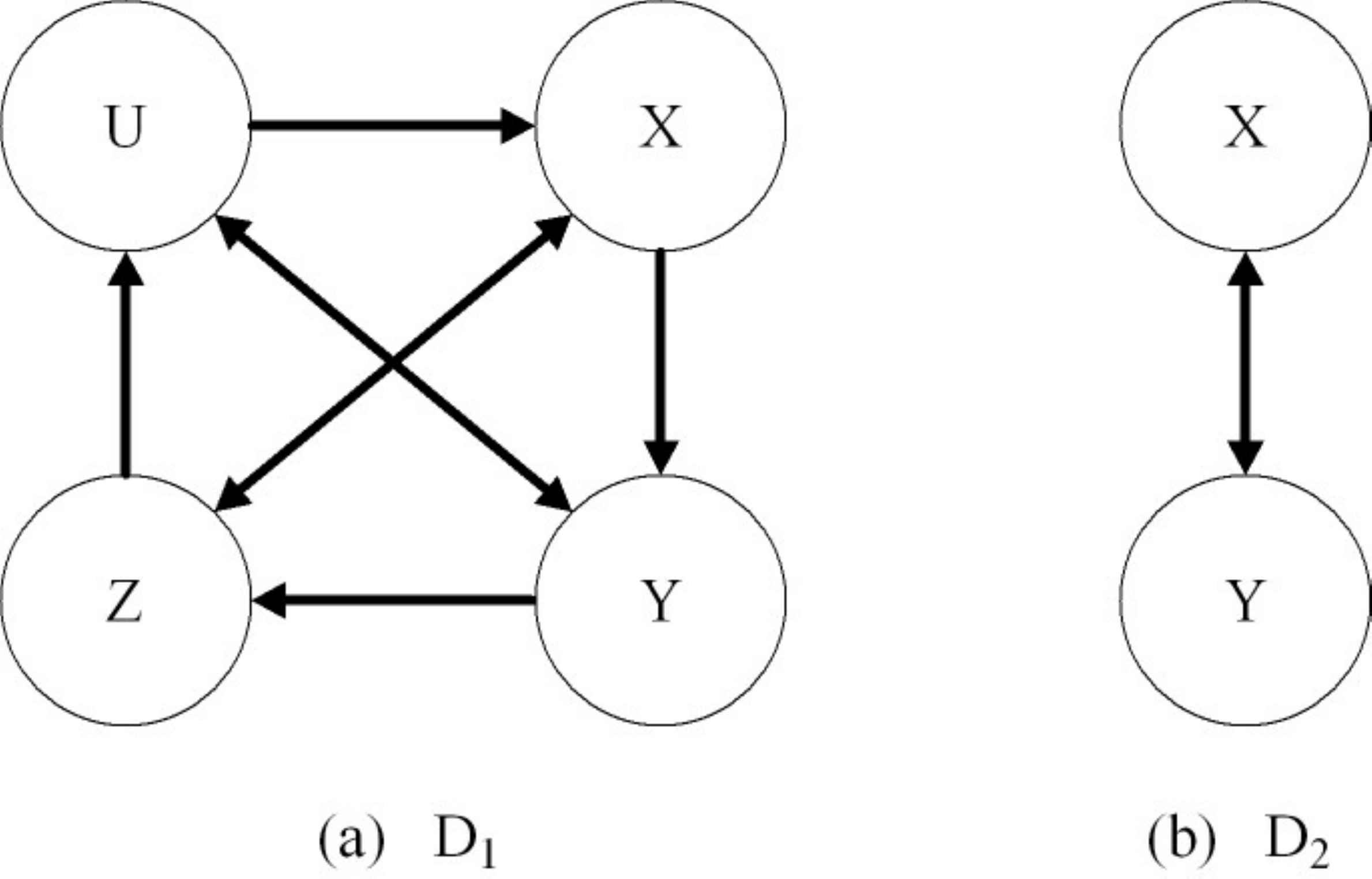}   
\caption{Counterexamples} 
\label{fig1}  
\end{figure}

Although the degree condition in Theorem \ref{thm1} is tight when $W=V(G)$ and $n_1=\cdots=n_k=2$, we believe that it can be improved if $n_i\geq 3$ for each $i$.  In fact, we  conjecture that  $2n/3$ is enough.

\begin{conj}\label{conj2}
Let $D$ be a digraph of order $n$ and $W\subseteq V(D)$. Suppose that $\delta^0(W)\geq 2n/3$ and $|W|=n_1+\cdots+n_k$ with $n_i\geq 3$ for each $i$. Then $D$ contains an arbitrary $W$-cycle-factor. That is, there are $k$ disjoint cycles $C_1,\ldots, C_k$ in $D$ such that $|V(C_i)\cap W|=n_i$ for all $i$.\end{conj}

The digraph $D_2$ in Fig.\ref{fig1} (b) shows that the minimum semi-degree in Conjecture \ref{conj2} is best possible.  $D_2$ consists two parts $X$ and $Y$, where $X=K^*_{2k-1}$ and $Y$ is an independent vertex set of size $k+1$. The bold arcs indicate complete domination in the direction shown. Clearly, $\delta^0(D_2)=2k-1=2n/3-1$, but it contains no $k$ disjoint cycles of length at least 3. Conjecture \ref{conj2} is supported by Theorems \ref{thmcyclable} and  \ref{thm5}. Moreover, in  \cite{Czygrinow2013}, Czygrinow, Kierstead and Molla  conjectured that every digraph $D$ of order $3k$ with $\delta^0(D)\geq2k$ contains $k$ disjoint triangles.  If Conjecture \ref{conj2} is true, then it implies this conjecture.


\begin{thebibliography}{99}

\bibitem{Abbasi1998} S. Abbasi, Ph. D Thesis (Rutgers University 1998).

\bibitem{Aigner1993}  M. Aigner, S. Brandt, Embedding arbitrary graphs of maximum degree two, J. Lond. Math. Soc. 1 (1993) 39-51.

\bibitem{Bang2009} J. Bang-Jensen, G. Gutin, Digraphs: Theory, Algorithms and Applications, 2nd edition, Springer-Verlag, London, 2009.

\bibitem{Berman1998} K. A. Berman, X. Liu, Cycles through large degree vertices in digraphs: A generalization of Meyniel's Theorem, J. Combin. Theory Ser. B 74 (1998) 20-27.

\bibitem{Bollobas1993} B. Bollob\'{a}s, G. Brightwell, Cycles through specified vertices, Combinatorica 13 (1993) 147-155.

\bibitem{Bondy1976} J. A. Bondy, U. S. R. Murty, Graph Theory, Springer London, 2008.

\bibitem{Chiba2018} S. Chiba, T. Yamashita, On directed 2-factors in digraphs and 2-factors containing perfect matchings in bipartite graphs, SIAM J. Discrete Math. 32 (2018) 394-409.

\bibitem{Corradi1963} K. Corr\'{a}di, A. Hajnal, On the maximal number of independent circuits in a graph, Acta. Math. Hung. 14 (1963) 423-439.

\bibitem{Czygrinow2013} A. Czygrinow, H. A. Kierstead, T. Molla, On directed versions of the Corr\'adi-Hajnal Corollary, Eur. J. Combin. 42 (2013) 1-14.

\bibitem{Dirac1952} G. A. Dirac, Some theorems on abstract graphs, P. Lond. Math. Soc. 2 (1952) 69-81.

\bibitem{Erdos1990} P. Erd\H{o}s, Some recent combinatorial problems, Technical Report, University of Bielefeld, Nov. 1990.

\bibitem{Fortune1980} S. Fortune, J. Hopcroft, J. Wyllie, The directed subgraph homeomorphism problem,   Theor. Comput. Sci. 10 (1980) 111-121.

\bibitem{Houri1962} A. Ghouila-Houri, Une condition suffisante d'existence d'un circuit hamiltonien,  C. R. Acad. Bulg. Sci.  25 (1960) 495-497. (In France)

\bibitem{Ore1960}  O. Ore, Note on hamilton circuits, Am. Math. Mon. 67 (1960) 55-55.


\bibitem{Sauer1978} N. Sauer and J. Spencer, Edge disjoint placement of graphs, J. Combin. Theory B,
25 (1978) 295-302.

\bibitem{Hong2000} H. Wang, Independent directed triangles in a directed graph, Graphs and Combin. 16 (2000) 453-462.

\bibitem{Hong2010} H. Wang, Proof of the Erd\H{o}s-Faudree conjecture on quadrilaterals, Graphs and Combin. 26 (2010) 833-877.

\bibitem{Hong2012} H. Wang, Disjoint 5-cycles in a graph, Discuss. Math. Graph Theory 32 (2012) 221-242.

\bibitem{Hong2013} H. Wang, Partition of a subset into two directed cycles, manuscript.

\bibitem{Hong2015} H. Wang, Partial degree conditions and cycle covering, J. Graph Therory 78 (2015) 267-304.

\bibitem{Zahar1984} M. H. El-Zahar, On circuits in graphs, Discrete Math. 50 (1984) 227-230.

\bibitem{Zhang2013} Z. B. Zhang, X. Zhang, X. Wen, Directed Hamilton cycles in digraphs and matching alternating Hamilton cycles in bipartite graphs, SIAM J. Discrete Math. 27 (2013) 274-289.
\end{thebibliography}
\end{document}